\documentclass[reqno]{amsart}%
\usepackage{float}
\usepackage{amsmath}
\usepackage{graphicx}
\usepackage{latexsym}
\usepackage{amsfonts}
\usepackage{amssymb}%
\theoremstyle{plain}
\newtheorem{theorem}{Theorem}

\newtheorem{lemma}[theorem]{Lemma}

\theoremstyle{definition}

\newtheorem{remark}[theorem]{Remark}
\newtheorem*{remark*}{Remark}

\begin{document}
\title{Local limit theorems for ladder epochs }
\author[Vatutin]{Vladimir A. Vatutin}
\address{Steklov Mathematical Institute RAS, Gubkin street 8, 19991, Moscow \\
Russia}
\email{vatutin@mi.ras.ru}
\author[Wachtel]{Vitali Wachtel}
\address{Weierstrass Institute for Applied Analysis and Stochastics, Mohrenstr.\ 39 \\
D--10117 Berlin, Germany}
\email{vakhtel@wias-berlin.de}
\urladdr{http://www.wias-berlin.de/$\sim$vakhtel}
\thanks{The first author was supported in part by the Russian Foundation for Basic
Research (grant 05-01-00035), INTAS (project 03-51-5018), and the program
"Contemporary problems of theoretical mathematics" of RAS}
\thanks{The second author was supported by the GIF}
 \subjclass{Primary 60\thinspace G\thinspace50; Secondary
60\thinspace G\thinspace40.}
\date{\today}

\begin{abstract}
Let $S_{0}=0,\{S_{n}\}_{n\geq1}$ be a random walk generated by a sequence of
i.i.d. random variables $X_{1},X_{2},...$ and let $\tau^{-}:=\min\left\{
n\geq1:\ S_{n}\leq0\right\}  $ and $\tau^{+}:=\min\left\{  n\geq
1:\ S_{n}>0\right\}  $. Assuming that the distribution of $X_{1}$ belongs to
the domain of attraction of an $\alpha$-stable law$,\alpha\neq1,$ we study the
asymptotic behavior of $\mathbb{P}(\tau^{\pm}=n)$ as $n\rightarrow\infty.$

\end{abstract}
\maketitle

\section{Introduction and main result}

Let $X,X_{1},X_{2},...$ be a sequence of independent identically distributed
random variables. Denote $S_{0}=0,S_{n}=X_{1}+X_{2}+...+X_{n}$. We assume
that
\[
\sum_{n=1}^{\infty}\frac{1}{n}\mathbb{P}(S_{n}>0)=\sum_{n=1}^{\infty}\frac
{1}{n}\mathbb{P}(S_{n}\leq0)=\infty.
\]
This condition means that $\{S_{n}\}_{n\geq0}$ is an oscillating random walk,
and, in particular, the stopping times
\[
\tau^{-}:=\min\left\{  n\geq1:\ S_{n}\leq0\right\}  \text{ \ and \ \ }\tau
^{+}:=\min\left\{  n\geq1:\ S_{n}>0\right\}
\]
are well-defined proper random variables. Furthermore, it follows from the
Wiener-Hopf factorization (see, for example, \cite[Theorem~8.9.1,~p.~376]%
{BGT}) that for all $z\in(0,1)$,
\begin{equation}
1-\mathbb{E}z^{\tau^{-}}=\exp\left\{  -\sum_{n=1}^{\infty}\frac{z^{n}}%
{n}\mathbb{P}(S_{n}\leq0)\right\}  \label{weak}%
\end{equation}
and
\begin{equation}
1-\mathbb{E}z^{\tau^{+}}=\exp\left\{  -\sum_{n=1}^{\infty}\frac{z^{n}}%
{n}\mathbb{P}(S_{n}>0)\right\}  . \label{strict}%
\end{equation}
Rogozin \cite{Rog71} has shown that the Spitzer condition
\begin{equation}
n^{-1}\sum_{k=1}^{n}\mathbb{P}\left(  S_{k}>0\right)  \rightarrow\rho
\in\left(  0,1\right)  \quad\text{as }n\rightarrow\infty\label{Spit}%
\end{equation}
holds if and only if $\tau^{+}$ belongs to the domain of attraction of a
spectrally positive stable law with parameter $\rho$. Since (\ref{weak}) and
(\ref{strict}) imply the equality
\[
(1-\mathbb{E}z^{\tau^{+}})(1-\mathbb{E}z^{\tau^{-}})=1-z\quad\text{for all
}z\in(0,1),
\]
one can deduce from the Rogozin result that (\ref{Spit}) holds if and only if
there exists a function $l(n)$ slowly varying at infinity such that, as
$n\rightarrow\infty$,
\begin{equation}
\mathbb{P}\left(  \tau^{-}>n\right)  \sim\frac{l(n)}{n^{1-\rho}}%
,\ \ \mathbb{P}\left(  \tau^{+}>n\right)  \sim\frac{1}{\Gamma(\rho
)\Gamma(1-\rho)n^{\rho}l(n)}. \label{integ1}%
\end{equation}
Doney \cite{Don95} proved that the Spitzer condition is equivalent to
\begin{equation}
\mathbb{P}\left(  S_{n}>0\right)  \rightarrow\rho\in\left(  0,1\right)
\quad\text{as }n\rightarrow\infty. \label{SpitDon}%
\end{equation}
Therefore, both relations in (\ref{integ1}) are valid under condition
(\ref{SpitDon}).

To get a more detailed information about the asymptotic properties of $l(x)$
it is necessary to impose additional hypotheses on the distribution of $X$.
Rogozin \cite{Rog71} has shown that $l(x)$ is asymptotically a constant if and
only if
\begin{equation}
\sum_{n=1}^{\infty}\frac{1}{n}\Bigl(\mathbb{P}\left(  S_{n}>0\right)
-\rho\Bigr)<\infty. \label{Rog}%
\end{equation}
It follows from the Spitzer-R\'{o}sen theorem (see
\cite[Theorem~8.9.23,~p.~382]{BGT}) that if $\mathbb{E}X^{2}$ is finite, then
(\ref{Rog}) holds with $\rho=1/2$, and, consequently,
\begin{equation}
\mathbb{P}(\tau^{\pm}>n)\sim\frac{C^{\pm}}{n^{1/2}}\quad\text{as }%
n\rightarrow\infty, \label{sym}%
\end{equation}
where $C^{\pm}$ are positive constants. If $\mathbb{E}X^{2}=\infty$ much less
is known about the form of $l(x)$. For instance, if the distribution of $X$ is
symmetric, then, clearly,
\begin{equation}
\left\vert \mathbb{P}\left(  S_{n}>0\right)  -\frac{1}{2}\right\vert =\frac
{1}{2}\mathbb{P}\left(  S_{n}=0\right)  . \label{simm}%
\end{equation}
Furthermore, according to \cite[Theorem III.9, p.~49]{Pet75}, there exists
$C>0$ such that for all $n\geq1$,
\[
\mathbb{P}\left(  S_{n}=0\right)  \leq\frac{C}{\sqrt{n}}.
\]
By this estimate and (\ref{simm}) we conclude that (\ref{Rog}) holds with
$\rho=1/2$. Thus, (\ref{sym}) is valid for all symmetric random walks.
Assuming that $\mathbb{P}(X>x)=\left(  x^{\alpha}l_{0}(x)\right)
^{-1},\ x>0,$ with $1<\alpha<2$ and $l_{0}(x)$ slowly varying at infinity,
Doney \cite{Don82a} established for a number of cases relationships between
the asymptotic behavior of $l_{0}(x)$ and $l(x)$ at infinity.

The aim of the present paper is to study the asymptotic behavior of the
probabilities $\mathbb{P}\left(  \tau^{\pm}=n\right)  $ as $n\rightarrow
\infty$.

We assume throughout that the distribution of $X$ is either non-lattice or
arithmetic with span $h>0$, i.e. the $h$ is the maximal number such that the
support of the distribution of $X$ is contained in the set $\left\{  kh,\text{
}k=0,\pm1,\pm2,...\right\}  .$

Let
\[
\mathcal{A}:=\{0<\alpha<1;\,|\beta|<1\}\cup\{1<\alpha<2;|\beta|\leq
1\}\cup\{\alpha=2,\beta=0\}
\]
be a subset in $\mathbb{R}^{2}.$ For $(\alpha,\beta)\in\mathcal{A}$ we write
$X\in\mathcal{D}\left(  \alpha,\beta\right)  $ if the distribution of $X$
belongs to the domain of attraction of a stable law with characteristic
function%
\begin{equation}
\Psi\mathbb{(}t\mathbb{)}:=\exp\left\{  -c|t|^{\,\alpha}\left(  1-i\beta
\frac{t}{|t|}\tan\frac{\pi\alpha}{2}\right)  \right\}  ,c>0, \label{std}%
\end{equation}
and, in addition, $\mathbb{E}X=0$ if $1<\alpha\leq2$. One can show
(see,$\ $\ for instance, \cite{Zol57}) that if $X\in\mathcal{D}\left(
\alpha,\beta\right)  $, then condition (\ref{SpitDon}) holds with
\begin{equation}
\rho=\frac{1}{2}+\frac{1}{\pi\alpha}\arctan\left(  \beta\tan\frac{\pi\alpha
}{2}\right)  \in(0,1). \label{ro}%
\end{equation}

Here is our main result.

\begin{theorem}
\label{TlocLad}Assume $X\in\mathcal{D}\left(  \alpha,\beta\right)  $. If
$\alpha\leq2$ and $\beta<1$, then, as $n\rightarrow\infty$,
\begin{equation}
\mathbb{P}\left(  \tau^{-}=n\right)  =(1-\rho)\frac{l(n)}{n^{2-\rho}}(1+o(1)).
\label{main}%
\end{equation}
In the case $\{1<\alpha<2,\,\beta=1\}$ equality $(\ref{main})$ remains valid
under the additional hypothesis
\begin{equation}
\int_{1}^{\infty}\frac{F(-x)}{x(1-F(x))}dx<\infty. \label{IntCond}%
\end{equation}

\end{theorem}

Denote $T^{-}:=\min\{n\geq1:\,S_{n}<0\}$ and set%
\begin{equation}
\Omega(z)=\exp\left\{  \sum_{n=1}^{\infty}\frac{z^{n}}{n}\mathbb{P}%
(S_{n}=0)\right\}  =:\sum_{k=0}^{\infty}\omega_{k}z^{k}. \label{defOmega}%
\end{equation}
The next statement relates the asymptotic behavior of $\mathbb{P}\left(
\tau^{-}=n\right)  $ and $\mathbb{P}\left(  T^{-}=n\right)  $.

\begin{theorem}
\label{propos} If $(\ref{main})$ holds, then
\[
\lim_{n\rightarrow\infty}\frac{\mathbb{P}\left(  T^{-}=n\right)  }%
{\mathbb{P}\left(  \tau^{-}=n\right)  }=\Omega(1).
\]

\end{theorem}

Applying Theorems \ref{TlocLad} and \ref{propos} to the random walk
$\{-S_{n}\}_{n\geq0}$, one can easily find an asymptotic representation for
$\mathbb{P}\left(  \tau^{+}=n\right)  $:

\begin{theorem}
\label{TlocLad'} Assume $X\in\mathcal{D}\left(  \alpha,\beta\right)  $. If
$\alpha\leq2$ and $\beta>-1$, then, as $n\rightarrow\infty$,
\begin{equation}
\mathbb{P}\left(  \tau^{+}=n\right)  =\frac{\rho}{\Gamma(\rho)\Gamma
(1-\rho)n^{1+\rho}l(n)}(1+o(1)). \label{main'}%
\end{equation}
In the case $\{1<\alpha<2,\,\beta=-1\}$ equality $(\ref{main'})$ remains valid
under the additional hypothesis
\begin{equation}
\label{IntCond'}\int_{1}^{\infty}\frac{1-F(x)}{xF(-x)}dx<\infty.
\end{equation}

\end{theorem}

In some special cases the asymptotic behavior of $\mathbb{P}\left(  \tau^{\pm
}=n\right)  $ as $n\rightarrow\infty$ is already known from the literature.
Eppel \cite{Epp} proved that if $\mathbb{E}X=0$ and $\mathbb{E}X^{2}$ is
finite, then
\begin{equation}
\mathbb{P}\left(  \tau^{\pm}=n\right)  \sim\frac{C^{\pm}}{n^{3/2}}.
\label{integ2}%
\end{equation}
Observe that in this case $\mathbb{E}X^{2}<\infty$ implies $X\in
\mathcal{D}(2,0)$.

Asymptotic representation (\ref{integ2}) is valid for all continuous symmetric
(implying $\rho=1/2$ in (\ref{SpitDon})) random walks (see \cite[Chapter XII,
Section 7]{FE}). Note that the restriction $X\in\mathcal{D}(\alpha,\beta)$ is
superfluous in this situation.

Recently Borovkov \cite{Bor04} has shown that if (\ref{Spit}) is valid and
\begin{equation}
n^{1-\rho}\Bigl(\mathbb{P}(S_{n}>0)-\rho\Bigr)\rightarrow const\in
(-\infty,\infty)\quad\text{as }n\rightarrow\infty, \label{BorCond}%
\end{equation}
then (\ref{main}) holds with $\ell(n)\equiv const\in(0,\infty)$. Proving the
mentioned result Borovkov does not assume that the distribution of $X$ is
taken from the domain of attraction of a stable law. However, he gives no
explanations how one can check the validity of (\ref{BorCond}) in the general situation.

Let $\chi^{+}:=S_{\tau^{+}}$ be the ascending ladder height. Alili and Doney
\cite[Remark~1,~p.~98]{AD99} have shown that (\ref{main'}) holds if
$\mathbb{E}\chi^{+}$ is finite. By Theorem 3 of \cite{Don82} the assumption
$\mathbb{E}\chi^{+}<\infty$ is equivalent to (\ref{IntCond'}), i.e. for the
case $\{1<\alpha<2,\beta=-1\}$ our Theorem~\ref{TlocLad'} is (implicitly)
contained in \cite{AD99} . Alili and Doney analyzed the distribution of
$\tau^{+}$ only. Clearly, one can easily derive the statement of our
Theorem~\ref{TlocLad} for the case $\{1<\alpha<2,\beta=1\}$ from their result
(for instance, applying Theorem~\ref{propos}). However, for these spectrally
one-sided cases we present an alternative proof, which clarifies the "typical"
behavior of the random walk on the events $\{\tau^{\pm}=n\}$. See Section 3.2
and Section 5 for more details.

%%%%%%%%%%%%%%%%%%%%%%%%%%%%%%%%%%%%%%%%%%%%%%%%%%%%%%%%%%%%%%%%%%%%%%%%%%%%%%%%%%%%%%%%%%
%%%%%%%%%%%%%%%%%%%%%%%%%%%%%%%%%%__PROOF__%%%%%%%%%%%%%%%%%%%%%%%%%%%%%%%%%%%%%%%%%%%%%%%
%%%%%%%%%%%%%%%%%%%%%%%%%%%%%%%%%%%%%%%%%%%%%%%%%%%%%%%%%%%%%%%%%%%%%%%%%%%%%%%%%%%%%%%%%%

\section{Auxiliary results}

\subsection{Notation}

In what follows we denote by $C,C_{1},C_{2},...$ finite positive constants
which may be \textit{different} from formula to formula and by $l(x),l_{1}%
(x),l_{2}(x)...$ functions slowly varying at infinity which are, as a rule,
\textit{fixed}.

For $x\geq0$ let%
\begin{align*}
\ B_{n}(x)  &  :=\mathbb{P}\left(  S_{n}\in(0,x];\tau^{-}>n\right)  ,\\
b_{n}(x)  &  :=B_{n}(x+1)-B_{n}(x)=\mathbb{P}\left(  S_{n}\in(x,x+1];\tau
^{-}>n\right)  .
\end{align*}
Introduce the renewal function
\[
H(x):=1+\sum_{k=1}^{\infty}\mathbb{P}\left(  \chi_{1}^{+}+...+\chi_{k}^{+}\leq
x\right)  ,\ x\geq0,\ H(x)=0,x<0,
\]
where $\left\{  \chi_{i}^{+}\right\}  _{i\geq1}$ is a sequence of i.i.d.
random variables distributed the same as $\chi^{+}$. Observe that by the
duality principle for random walks for$\ \ x\geq0$
\begin{align}
1+\sum_{j=1}^{\infty}B_{j}(x)  &  =1+\sum_{j=1}^{\infty}\mathbb{P}\left(
S_{j}\in(0,x];\tau^{-}>j\right) \nonumber\\
&  =1+\sum_{j=1}^{\infty}\mathbb{P}\left(  S_{j}\in(0,x];S_{j}>S_{0}%
,S_{j}>S_{1},...,S_{j}>S_{j-1}\right) \nonumber\\
&  =H(x). \label{Dual}%
\end{align}

In the sequel we deal rather often with slowly varying functions and,
following Doney \cite{Don82}, say that a slowly varying function $l^{\ast}(x)$
is an $\alpha$-conjugate of a slowly varying function $l^{\ast\ast}(x)$ when
the following relations are valid
\[
y\sim x^{\alpha}l^{\ast}(x)\text{ as }x\rightarrow\infty\text{ if and only if
}x\sim y^{1/\alpha}l^{\ast\ast}(y).
\]

It is known that if $X\in\mathcal{D}\left(  \alpha,\beta\right)  $ with
$\alpha\in(0,2),$ and $F(x):=\mathbb{P}\left(  X\leq x\right)  $, then%
\begin{equation}
1-F(x)+F(-x)\sim\frac{1}{x^{\alpha}l_{0}(x)}\quad\text{as }x\rightarrow\infty,
\label{Tailtwo}%
\end{equation}
where $l_{0}(x)$ is a function slowly varying at infinity. Besides, for
$\alpha\in(0,2)$,%
\begin{equation}
\frac{F(-x)}{1-F(x)+F(-x)}\rightarrow q,\quad\frac{1-F(x)}{1-F(x)+F(-x)}%
\rightarrow p\quad\text{as }x\rightarrow\infty, \label{tailF}%
\end{equation}
with $p+q=1$ and $\beta=p-q$ in (\ref{std}). Let $\left\{  c_{n}\right\}
_{n\geq1}$ be a sequence specified by the relation%
\begin{equation}
c_{n}:=\inf\left\{  x\geq0:1-F(x)+F(-x)\leq n^{-1}\right\}  . \label{Defa}%
\end{equation}
In view of (\ref{Tailtwo}) this sequence is regularly varying at infinity with
index $\alpha^{-1}$, i.e.
\begin{equation}
c_{n}=n^{1/\alpha}l_{1}(n), \label{asyma}%
\end{equation}
where $l_{1}(x)$ is a slowly varying function being an $\alpha$-conjugate of
$l_{0}(x)$:
\begin{equation}
c_{n}^{\alpha}l_{0}(c_{n})\sim n\quad\text{as }n\rightarrow\infty.
\label{Repa}%
\end{equation}
Moreover,%
\[
\frac{S_{n}}{c_{n}}\overset{d}{\rightarrow}Y_{\alpha}\quad\text{as
}n\rightarrow\infty,
\]
where $Y_{\alpha}$ is a random variable obeying an $\alpha-$stable law.

For the case $\alpha=2$ the normalizing sequence $\left\{  c_{n}\right\}
_{n\geq1}$ requires a special description. Let $V(x)=\int_{-x}^{x}y^{2}dF(x)$
be the truncated variance of $X$. Clearly, $\liminf_{x\rightarrow\infty
}V(x)>0$ for every nondegenerate random variable $X$. Furthermore, it is known
(\cite{FE}, Chapter XVII, Section 5) that $X\in\mathcal{D}(2,0)$ if and only
if $V(x)$ varies slowly at infinity. In this case the normalizing sequence
$c_{n}$ satisfies
\begin{equation}
\frac{V(c_{n})}{c_{n}^{2}}\sim\frac{C}{n}\quad\text{as }n\rightarrow\infty.
\label{102}%
\end{equation}
The last relation means that (\ref{asyma}) holds with $\alpha=2$ and
$l_{1}(x)$ is a $2$-conjugate of $1/V(x)$. Besides, \
\begin{equation}
\lim_{x\rightarrow\infty}\frac{x^{2}(1-F(x)+F(-x))}{V(x)}=0. \label{ProV}%
\end{equation}
%%%%%%%%%%%%%%%%%%%%%%%%%%%%%%%%%%%%%%%%%%%%%%%%%%%%%%%%%%%%%%%%%%%%%%%%%%%%%

\subsection{Basic lemmas}

Now we formulate a number of results concerning the distributions of the
random variables $\tau^{-},\tau^{+}$ and $\chi^{+}$. $\ $\ Recall that a
random variable $\zeta$ is called relatively stable if there exists a
nonrandom sequence $d_{n}\rightarrow\infty$ as $n\rightarrow\infty$ such that
\[
\frac{1}{d_{n}}\sum_{k=1}^{n}\zeta_{k}\overset{p}{\rightarrow}1\text{ as
}n\rightarrow\infty,
\]
where $\zeta_{k}\overset{d}{=}\zeta,\ k=1,2,...$ and are independent.

\begin{lemma}
\label{Lrenewal}$\mathrm{(}$see \cite{Rog71} and \cite[Theorem~9]%
{Don85}$\mathrm{)}$ Assume $X\in\mathcal{D}(\alpha,\beta)$. Then, as
$x\rightarrow\infty$,%
\begin{equation}
\mathbb{P}\left(  \chi^{+}>x\right)  \sim\frac{1}{x^{\alpha\rho}l_{2}%
(x)}\ \text{if }\alpha\rho<1, \label{asymXi}%
\end{equation}
and $\chi^{+}$ is relatively stable if $\alpha\rho=1.$
\end{lemma}

\begin{lemma}
\label{Renew2} Suppose $X\in\mathcal{D}(\alpha,\beta)$. If $\alpha\rho<1$,
then, as $x\rightarrow\infty$,%
\begin{equation}
H(x)\sim\frac{x^{\alpha\rho}l_{2}(x)}{\Gamma(1-\alpha\rho)\Gamma(1+\alpha
\rho)}. \label{RenStand}%
\end{equation}
If $\alpha\rho=1$, then, as $x\rightarrow\infty$,
\begin{equation}
H(x)\sim xl_{3}(x), \label{RenewRelat}%
\end{equation}
where%
\[
l_{3}(x):=\left(  \int_{0}^{x}\mathbb{P}\left(  \chi^{+}>y\right)  dy\right)
^{-1},\text{ \ }x>0.
\]
In addition, there exists a constant $C>0$ such that in both cases
\begin{equation}
H(c_{n})\leq Cn^{\rho}l(n)\quad\text{for all }n\geq1. \label{H-bound}%
\end{equation}

\end{lemma}

\begin{proof}
If $\alpha\rho<1$, then by \cite[Chapter~XIV,~formula~(3.4)]{FE}
\[
H(x)\sim\frac{1}{\Gamma(1-\alpha\rho)\Gamma(1+\alpha\rho)}\frac{1}%
{\mathbb{P}(\chi^{+}>x)}\quad\text{as }x\rightarrow\infty.
\]
Hence, recalling (\ref{asymXi}), we obtain (\ref{RenStand}).

If $\alpha\rho=1$, then (\ref{RenewRelat}) follows from Theorem 2 in
\cite{Rog71}.

Let us demonstrate the validity of \ (\ref{H-bound}). We know from
\cite{Rog71} (see also \cite{GOT82}) that $\tau^{+}\in\mathcal{D}(\rho,1)$
under the conditions of the lemma and, in addition, $\chi^{+}\in
\mathcal{D}(\alpha\rho,1)$ if $\alpha\rho<1$. This means, in particular, that
for sequences $\{a_{n}\}_{n\geq1}$ and $\{b_{n}\}_{n\geq1}$ specified by
\begin{equation}
\mathbb{P}(\tau^{+}>a_{n})\sim\frac{1}{n}\quad\text{and}\quad\mathbb{P}%
(\chi^{+}>b_{n})\sim\frac{1}{n}\quad\text{as }n\rightarrow\infty,
\label{an-asym}%
\end{equation}
and vectors $\{(\tau_{k}^{+},\chi_{k}^{+})\}_{k\geq1},$ being independent
copies of $(\tau^{+},\chi^{+}),$ we have
\begin{equation}
\frac{1}{a_{n}}\sum_{k=1}^{n}\tau_{k}^{+}\overset{d}{\rightarrow}Y_{\rho}%
\quad\text{and}\quad\frac{1}{b_{n}}\sum_{k=1}^{n}\chi_{k}^{+}\overset
{d}{\rightarrow}Y_{\alpha\rho}\text{ \ \ as \ }n\rightarrow\infty. \label{**}%
\end{equation}
Moreover, it was established by Doney (see Lemma in \cite{Don85}, p. 358)
that
\begin{equation}
b_{n}\sim Cc_{[a_{n}]}\text{ as }n\rightarrow\infty, \label{doney}%
\end{equation}
where $[x]$ stands for the integer part of $x$. Therefore, $c_{n}\sim
Cb_{[a^{-1}(n)]}$, where, with a slight abuse of notation, $a^{-1}(n)$ is the
inverse function to $a_{n}$. Hence, on account of (\ref{an-asym}),
\begin{align}
\mathbb{P}(\chi^{+}>c_{n})  &  \sim C_{1}\mathbb{P}(\chi^{+}>b_{[a^{-1}%
(n)]})\sim\frac{C_{1}}{a^{-1}(n)}\nonumber\label{***}\\
&  \sim C_{2}\mathbb{P}(\tau^{+}>a_{[a^{-1}(n)]})\sim C_{3}\mathbb{P}(\tau
^{+}>n)\sim\frac{C_{4}}{n^{\rho}l(n)}.
\end{align}
This proves (\ref{H-bound}) for $\alpha\rho<1$.

If $\alpha\rho=1$, then, instead of the second equivalence in (\ref{an-asym}),
one should define $b_{n}$ by
\[
\frac{1}{b_{n}}\int_{0}^{b_{n}}\mathbb{P}(\chi^{+}>y)dy\sim\frac{1}{n}%
\quad\text{as }n\rightarrow\infty
\]
(see \cite[p.~595]{Rog71}). In this case the second convergence in (\ref{**})
transforms to
\[
\frac{1}{b_{n}}\sum_{k=1}^{n}\chi_{k}^{+}\overset{p}{\rightarrow}%
1\quad\text{as }n\rightarrow\infty,
\]
while (\ref{***}) should be changed to
\begin{align*}
\frac{1}{c_{n}}\int_{0}^{c_{n}}\mathbb{P}(\chi^{+}>y)dy  &  \sim\frac{C_{1}%
}{b_{[a^{-1}(n)]}}\int_{0}^{b_{[a^{-1}(n)]}}\mathbb{P}(\chi^{+}>y)dy\sim
\frac{C_{1}}{a^{-1}(n)}\\
&  \sim C_{1}\mathbb{P}(\tau^{+}>a_{[a^{-1}(n)]})\sim C_{2}\mathbb{P}(\tau
^{+}>n)\sim\frac{C_{3}}{n^{\rho}l(n)}.
\end{align*}
The lemma is proved.
\end{proof}

The next result is a part of Corollary 3 in \cite{Don82}.

\begin{lemma}
\label{beta1} Assume $X\in\mathcal{D}(\alpha,1)$ with $1<\alpha<2$ $($implying
$\rho=1-\alpha^{-1})$. Then
\[
\mathbb{P}(\tau^{-}>n)\sim\frac{C}{c_{n}}\sim\frac{C}{n^{1/\alpha}l_{1}%
(n)}\ \ \ \text{\ as }n\rightarrow\infty
\]
if and only if
\[
\int_{1}^{\infty}\frac{F(-x)}{x(1-F(x))}dx<\infty.
\]

\end{lemma}

Now we prove a useful result which may be viewed as a statement concerning
"small" deviations of $S_{n}$ on the set $\left\{  \tau^{-}>n\right\}  $. \

Let $h$ be the span and $g_{\alpha,\beta}(x)$ be the density of a stable
distribution with parameters $\alpha$ and $\beta$ in (\ref{std}) (we agree to
consider $h=0$ for non-lattice distributions). For a set $A$ taken from the
Borel $\sigma$-algebra on $(0,\infty)$ denote
\[
\mu(A)=g_{\alpha,\beta}(0)\int_{A}H(x-h)\nu(dx),
\]
where $\nu$ is the counting measure on $\{h,2h,3h,\ldots\}$ in the arithmetic
case and the Lebesgue measure on $(0,\infty)$ in the non-lattice case.

\begin{lemma}
\label{L_Conv} Suppose $X\in\mathcal{D}(\alpha,\beta)$. Then
\begin{equation}
\lim_{n\rightarrow\infty}nc_{n}\mathbb{P}(S_{n}\in A;\,\tau^{-}>n)=\mu(A)
\label{Conv}%
\end{equation}
for any $A$ taken from the Borel $\sigma$-algebra on $(0,\infty)$.
\end{lemma}

\begin{proof}
Assume first that the distribution of $X$ is non-lattice. Using the Stone
local limit theorem (see, for instance, \cite[Section~8.4,~p.~351]{BGT}) it is
not difficult to show that for $\lambda>0$,
\begin{equation}
\lim_{n\rightarrow\infty}c_{n}\mathbb{E}\left(  e^{-\lambda S_{n}}%
;\,S_{n}>0\right)  =g_{\alpha,\beta}(0)\int_{0}^{\infty}e^{-\lambda y}%
dy=\frac{g_{\alpha,\beta}(0)}{\lambda}. \label{sm1}%
\end{equation}
Set%
\begin{equation}
G(\lambda):=\sum_{n=1}^{\infty}\frac{\mathbb{E}\left(  e^{-\lambda S_{n}%
};\,S_{n}>0\right)  }{n} \label{sm2}%
\end{equation}
and specify a sequence of measures
\[
\mu_{n}(dx):=nc_{n}\mathbb{P(}S_{n}\in dx;\,\tau^{-}>n),\quad n\geq1.
\]
Since $\left\{  c_{n}\right\}  _{n\geq1}$ varies regularly and (\ref{sm1}) is
valid, applying Theorem 2 from \cite{CNW} to the equality
\begin{equation}
\sum_{n=0}^{\infty}z^{n}\mathbb{E}\left(  e^{-\lambda S_{n}};\,\tau
^{-}>n\right)  =\exp\left\{  \sum_{n=1}^{\infty}\frac{z^{n}}{n}\mathbb{E}%
\left(  e^{-\lambda S_{n}};\,S_{n}>0\right)  \right\}  \label{RR1}%
\end{equation}
shows that for all $\lambda>0$,
\begin{align}
\lim_{n\rightarrow\infty}nc_{n}\mathbb{E}\left(  e^{-\lambda S_{n}};\,\tau
^{-}>n\right)   &  =\lim_{n\rightarrow\infty}\int_{0}^{\infty}e^{-\lambda
x}\mu_{n}(dx)\nonumber\\
&  =\frac{g_{\alpha,\beta}(0)}{\lambda}\exp\left\{  G(\lambda)\right\}  .
\label{sm3}%
\end{align}
It follows from (\ref{RR1}) that
\begin{align*}
\frac{g_{\alpha,\beta}(0)}{\lambda}\exp\left\{  G(\lambda)\right\}   &
=\frac{g_{\alpha,\beta}(0)}{\lambda}\Bigl(1+\sum_{k=1}^{\infty}\mathbb{E}%
\left(  e^{-\lambda S_{k}};\,\tau^{-}>k\right)  \Bigr)\\
&  =\frac{g_{\alpha,\beta}(0)}{\lambda}+\frac{g_{\alpha,\beta}(0)}{\lambda
}\int_{0}^{\infty}e^{-\lambda x}\left(  \sum_{k=1}^{\infty}\mathbb{P}\left(
S_{k}\in dx;\,\tau^{-}>k\right)  \right) \\
&  =\frac{g_{\alpha,\beta}(0)}{\lambda}+\frac{g_{\alpha,\beta}(0)}{\lambda
}\int_{0}^{\infty}e^{-\lambda x}\left(  \sum_{j=1}^{\infty}\mathbb{P}(\chi
_{1}^{+}+\ldots+\chi_{j}^{+}\in dx)\right)  ,
\end{align*}
where at the last step we have used the duality principle. Integrating by
parts and recalling the definition of $H(x)$, we get
\begin{align}
\frac{g_{\alpha,\beta}(0)}{\lambda}\exp\left\{  G(\lambda)\right\}   &
=\frac{g_{\alpha,\beta}(0)}{\lambda}+g_{\alpha,\beta}(0)\int_{0}^{\infty
}e^{-\lambda x}(H(x)-1)dx\nonumber\\
&  =g_{\alpha,\beta}(0)\int_{0}^{\infty}e^{-\lambda x}H(x)dx. \label{sm4}%
\end{align}
Combining (\ref{sm3}) and (\ref{sm4}) and using the continuity theorem for
Laplace transforms, we obtain (\ref{Conv}) for non-lattice distributions.

In the arithmetic case we have by the Gnedenko local limit theorem%
\begin{equation}
\label{sm5}\lim_{n\rightarrow\infty}c_{n}\mathbb{E}\left(  e^{-\lambda S_{n}%
};\,S_{n}>0\right)  =g_{\alpha,\beta}(0)\sum_{k=1}^{\infty}e^{-\lambda
hk}=\frac{g_{\alpha,\beta}(0)e^{-\lambda h}}{1-e^{-\lambda h}}.
\end{equation}
Proceeding as by the derivation of (\ref{sm4}), we obtain
\begin{align*}
\frac{g_{\alpha,\beta}(0)e^{-\lambda h}}{1-e^{-\lambda h}}\exp\left\{
G(\lambda)\right\}   &  =\frac{g_{\alpha,\beta}(0)e^{-\lambda h}%
}{1-e^{-\lambda h}}\Bigl(1+\sum_{k=1}^{\infty}\mathbb{E}\left(  e^{-\lambda
S_{k}};\,\tau^{-}>k\right)  \Bigr)\\
&  \hspace{-1cm} =\frac{g_{\alpha,\beta}(0)e^{-\lambda h}}{1-e^{-\lambda h}%
}+\frac{g_{\alpha,\beta}(0)e^{-\lambda h}}{1-e^{-\lambda h}}\sum_{j=1}%
^{\infty}e^{-\lambda hj}\left(  H(hj)-H(hj-h)\right) \\
&  \hspace{-1cm} =g_{\alpha,\beta}(0)e^{-\lambda h}\sum_{j=0}^{\infty
}e^{-\lambda hj}H(hj)=g_{\alpha,\beta}(0)\sum_{k=1}^{\infty}e^{-\lambda
hk}H(hk-h).
\end{align*}
This, together with (\ref{sm5}), finishes the proof of the lemma.
\end{proof}

\begin{lemma}
\label{L0}Under the conditions of Theorem \textrm{\ref{TlocLad}} for any
$\alpha\in(0,2)$ there exists $C>0$ such that for all $y>0$ and all $n\geq1$,%
\begin{equation}
b_{n}(y)\leq\frac{C}{c_{n}}\frac{l(n)}{n^{1-\rho}} \label{uupb}%
\end{equation}
and%
\begin{equation}
B_{n}(y)\leq\frac{C\left(  y+1\right)  }{c_{n}}\frac{l(n)}{n^{1-\rho}}.
\label{UpB}%
\end{equation}

\end{lemma}

\begin{proof}
For $n=1$ the statement of the lemma is obvious. Let $\left\{  S_{n}^{\ast
}\right\}  _{n\geq0}$ be a random walk distributed as $\left\{  S_{n}\right\}
_{n\geq0}$ and independent of it. One can easily check that for each $n\geq2$,%
\begin{align}
b_{n}(y)  &  =\mathbb{P}\left(  y<S_{n}\leq y+1;\tau^{-}>n\right) \nonumber\\
&  =\int_{0}^{\infty}\mathbb{P}\Bigl(y-S_{[n/2]}<S_{n}-S_{[n/2]}\leq
y+1-S_{[n/2]};S_{[n/2]}\in dz;\tau^{-}>n\Bigr)\nonumber\\
&  \leq\int_{0}^{\infty}\mathbb{P}\Bigl(y-z<S_{n-[n/2]}^{\ast}\leq
y+1-z;S_{[n/2]}\in dz;\tau^{-}>[n/2]\Bigr)\nonumber\\
&  \leq\mathbb{P}\Bigl(\tau^{-}>[n/2]\Bigr)\sup_{z}\mathbb{P}%
\Bigl(z<S_{n-[n/2]}^{\ast}\leq z+1\Bigr). \label{Sm1}%
\end{align}
Since the density of any $\alpha$-stable law is bounded, it follows from the
Gnedenko and Stone local limit theorems that if the distribution of $X$ is
either arithmetic or non-lattice, then there exists a constant $C>0$ such that
for all $n\geq1$ and all $z\geq0$,
\begin{equation}
\mathbb{P}\left(  S_{n}\in(z,z+\Delta]\right)  \leq\frac{C\Delta}{c_{n}}.
\label{EstS1}%
\end{equation}
Hence it follows, in particular, that, for any $z>0$,
\begin{equation}
\mathbb{P}\left(  S_{n}\in(0,z]\right)  \leq\frac{C(z+1)}{c_{n}}.
\label{EstS2}%
\end{equation}
Substituting (\ref{EstS1}) into (\ref{Sm1}), and recalling (\ref{asyma}) and
properties of regularly varying functions, we get (\ref{uupb}). Estimate
(\ref{UpB}) follows from (\ref{uupb}) by summation.
\end{proof}

\begin{lemma}
\label{L1}Under the conditions of Theorem \ref{TlocLad} for any $\alpha
\in(1,2]$ there exists $C>0$ such that for all $n\geq1$ and all $x>0$,%
\begin{equation}
b_{n}(x)\leq C\left(  \frac{H(x+1)}{nc_{n}}+\frac{l(n)}{n^{1-\rho}}\frac
{x+1}{c_{n}^{2}}\right)  \label{bsmall2}%
\end{equation}
and
\begin{equation}
B_{n}(x)\leq C\left(  \frac{\left(  x+1\right)  H(x+1)}{nc_{n}}+\frac
{l(n)}{n^{1-\rho}}\frac{\left(  x+1\right)  ^{2}}{c_{n}^{2}}\right)  .
\label{Bbig2}%
\end{equation}

\end{lemma}

\begin{proof}
According to formula (5) in \cite{Epp},
\begin{equation}
nB_{n}(x)=\mathbb{P}\left(  S_{n}\in(0,x]\right)  +\sum_{k=1}^{n-1}\int
_{0}^{x}B_{n-k}(x-y)\mathbb{P}\left(  S_{k}\in dy\right)  . \label{Eppident}%
\end{equation}
Hence we get%
\begin{align}
nb_{n}(x)  &  =\mathbb{P}\left(  S_{n}\in(x,x+1]\right)  +\sum_{k=1}^{n-1}%
\int_{0}^{x}b_{n-k}(x-y)\mathbb{P}\left(  S_{k}\in dy\right) \nonumber\\
&  \qquad+\sum_{k=1}^{n-1}\int_{x}^{x+1}B_{n-k}(x+1-y)\mathbb{P}\left(
S_{k}\in dy\right)  . \label{Smallb}%
\end{align}
Using (\ref{uupb}), (\ref{EstS2}), (\ref{asyma}), the inequality $1/\alpha<1$
and properties of slowly varying functions, we deduce
\begin{align}
\sum_{k=1}^{[n/2]}\int_{0}^{x}b_{n-k}(x-y)\mathbb{P}\left(  S_{k}\in
dy\right)   &  \leq C\sum_{k=1}^{[n/2]}\frac{l(n-k)}{c_{n-k}\left(
n-k\right)  ^{1-\rho}}\mathbb{P}\left(  S_{k}\in\lbrack0,x]\right) \nonumber\\
&  \leq C_{1}\left(  x+1\right)  \sum_{k=1}^{[n/2]}\frac{1}{c_{k}}%
\frac{l(n-k)}{c_{n-k}\left(  n-k\right)  ^{1-\rho}}\nonumber\\
&  \leq C_{2}\frac{x+1}{c_{n}}\frac{l(n)}{n^{1-\rho}}\sum_{k=1}^{[n/2]}%
\frac{1}{c_{k}}\nonumber\\
&  \leq C_{3}\left(  x+1\right)  \frac{n^{\rho}l(n)}{c_{n}^{2}}.
\label{EstSmall}%
\end{align}
On the other hand, in view of (\ref{EstS1}) and monotonicity of $B_{k}(x)$ in
$x$ we conclude (assuming that $x$ is integer without loss of generality and
letting $B_{k}(-1)=0$ and $H(-1)=0$) that
\begin{align*}
&  \sum_{k=[n/2]+1}^{n}\int_{0}^{x}b_{n-k}(x-y)\mathbb{P}\left(  S_{k}\in
dy\right) \\
&  \leq\sum_{k=[n/2]+1}^{n}\sum_{j=0}^{x}\left(  B_{n-k}(x-j+1)-B_{n-k}%
(x-j-1)\right)  \mathbb{P}\left(  S_{k}\in(j,j+1]\right) \\
&  \leq\sum_{k=[n/2]+1}^{n}\sum_{j=0}^{x}\left(  B_{n-k}(x-j+1)-B_{n-k}%
(x-j-1)\right)  \frac{C}{c_{k}}\\
&  \leq\frac{C}{c_{n}}\sum_{j=0}^{x}\sum_{k=0}^{\infty}\left(  B_{k}%
(x-j+1)-B_{k}(x-j-1)\right) \\
&  =\frac{C}{c_{n}}\sum_{j=0}^{x}\left(  H(x-j+1)-H(x-j-1)\right) \\
&  \leq\frac{C}{c_{n}}\left(  H(x)+H(x+1)\right)  \leq\frac{2C}{c_{n}}H(x+1),
\end{align*}
where for the intermediate equality we have used (\ref{Dual}). This gives
\begin{equation}
\sum_{k=[n/2]+1}^{n}\int_{0}^{x}b_{n-k}(x-y)\mathbb{P}\left(  S_{k}\in
dy\right)  \leq\frac{C}{c_{n}}H(x+1). \label{EstBig}%
\end{equation}
Since $x\mapsto B_{n}(x)$ increases for every $n$,
\begin{equation}
\sum_{k=1}^{n-1}\int_{x}^{x+1}B_{n-k}(x+1-y)\mathbb{P}\left(  S_{k}\in
dy\right)  \leq\sum_{k=1}^{n-1}B_{n-k}(1)\mathbb{P}(S_{k}\in(x,x+1]).
\label{added1}%
\end{equation}
Further, in view of (\ref{UpB}) and (\ref{EstS1}) we have
\begin{equation}
\sum_{k=1}^{[n/2]}B_{n-k}(1)\mathbb{P}(S_{k}\in(x,x+1])\leq\frac{C_{1}}{c_{n}%
}\frac{l(n)}{n^{1-\rho}}\sum_{k=1}^{[n/2]}\frac{1}{c_{k}}\leq\frac
{C_{2}n^{\rho}l(n)}{c_{n}^{2}}. \label{added2}%
\end{equation}
Using (\ref{EstS1}) once again yields
\begin{equation}
\sum_{k=[n/2]+1}^{n-1}B_{n-k}(1)\mathbb{P}(S_{k}\in(x,x+1])\leq\frac{C}{c_{n}%
}\sum_{k=[n/2]+1}^{n-1}B_{n-k}(1)\leq\frac{C}{c_{n}}H(1). \label{added3}%
\end{equation}
Substituting (\ref{added2}) and (\ref{added3}) into the right hand side of
(\ref{added1}), we obtain the upper bound
\begin{equation}
\sum_{k=1}^{n-1}\int_{x}^{x+1}B_{n-k}(x+1-y)\mathbb{P}\left(  S_{k}\in
dy\right)  \leq C\Bigl(\frac{n^{\rho}l(n)}{c_{n}^{2}}+\frac{1}{c_{n}}\Bigr).
\label{EstimRem}%
\end{equation}
Combining (\ref{EstSmall}), (\ref{EstBig}), (\ref{EstimRem}), (\ref{EstS1})
and (\ref{Smallb}) proves (\ref{bsmall2}). Observing that $H(x)$ is
nondecreasing and integrating (\ref{bsmall2}), we get estimate (\ref{Bbig2}).
\end{proof}

%%%%%%%%%%%%%%%%%%%%%%%%%%%%%%%%%%%%%%%%%%%%%%%%%%%%%%%%%%%%%%%%%%%%%%%%%%%%%%%%%%

To prove Theorem \ref{TlocLad} in the case $\alpha=2$ we need the following
technical lemma which may be known from the literature.

\begin{lemma}
\label{Ltec}Let $w(n)$ be a monotone increasing function. If, \ for some
$\gamma>0,$ there exist slowly varying functions $l^{\ast}(n)$ and
$l^{\ast\ast}(n)$ such that, as $n\rightarrow\infty$,%
\[
\sum_{k=n}^{\infty}\frac{w(k)}{k^{\gamma+1}l^{\ast}(k)}\sim\frac{1}{n^{\gamma
}l^{\ast\ast}(n)},
\]
then, as $n\rightarrow\infty$,%
\[
w(n)\sim\gamma\frac{l^{\ast}(n)}{l^{\ast\ast}(n)}.
\]

\end{lemma}

\begin{proof}
Let, for this lemma only, $r_{i}(n),n=1,2,...;i=1,2,3,4$ be sequences of real
numbers vanishing as $n\rightarrow\infty.$ For $\Delta\in(0,1)$ we have by
monotonicity of $w(n)$ and properties of slowly varying functions%
\begin{align*}
w(\left[  \Delta n\right]  )\sum_{k=\left[  \Delta n\right]  }^{n}\frac
{1}{k^{\gamma+1}l^{\ast}(k)}  &  =w(\left[  \Delta n\right]  )\frac
{1+r_{2}(n)}{\gamma n^{\gamma}l^{\ast}(n)}\left(  \Delta^{-\gamma}-1\right) \\
&  \leq\sum_{k=\left[  \Delta n\right]  }^{n}\frac{w(k)}{k^{\gamma+1}l^{\ast
}(k)}=\frac{1+r_{1}(n)}{n^{\gamma}l^{\ast\ast}(n)}\left(  \Delta^{-\gamma
}-1\right) \\
&  \leq w(n)\sum_{k=\left[  \Delta n\right]  }^{n}\frac{1}{k^{\gamma+1}%
l^{\ast}(k)}\\
&  =w(n)\frac{1+r_{2}(n)}{\gamma n^{\gamma}l^{\ast}(n)}\left(  \Delta
^{-\gamma}-1\right)  .
\end{align*}
Hence it follows that%
\[
w(\left[  \Delta n\right]  )\leq\frac{1+r_{1}(n)}{1+r_{2}(n)}\frac{\gamma
l^{\ast}(n)}{l^{\ast\ast}(n)}\leq w(n)
\]
and, therefore,%
\[
\frac{1+r_{1}(n)}{1+r_{2}(n)}\frac{\gamma l^{\ast}(n)}{l^{\ast\ast}(n)}\leq
w(n)\leq\frac{1+r_{3}(\left[  n\Delta^{-1}\right]  )}{1+r_{4}(\left[
n\Delta^{-1}\right]  )}\frac{\gamma l^{\ast}(\left[  n\Delta^{-1}\right]
)}{l^{\ast\ast}(\left[  n\Delta^{-1}\right]  )}.
\]
Since $l^{\ast}$ and $l^{\ast\ast}$ are slowly varying functions, we get%
\[
\lim_{n\rightarrow\infty}\frac{w(n)l^{\ast\ast}(n)}{\gamma l^{\ast}(n)}=1,
\]
as desired.
\end{proof}

\begin{remark}
\label{Rem22}By the same arguments one can show that if $w(x)$ is a monotone
increasing function and,\ for some $\gamma>0,$ there exist slowly varying
functions $l^{\ast}(x)$ and $l^{\ast\ast}(x)$ such that, as $x\rightarrow
\infty$,%
\[
\int_{x}^{\infty}\frac{w(y)dy}{y^{\gamma+1}l^{\ast}(y)}\sim\frac{1}{x^{\gamma
}l^{\ast\ast}(x)},
\]
then, as $x\rightarrow\infty$,%
\[
w(x)\sim\gamma\frac{l^{\ast}(x)}{l^{\ast\ast}(x)}.
\]

\end{remark}

\section{Proof of Theorem \ref{TlocLad}}

\subsection{Proof of Theorem \ref{TlocLad} for $\{0<\alpha<2,\,\beta
<1\}\cap\left\{  \alpha\neq1\right\}  $}

For a fixed $\varepsilon\in(0,1)$ write
\begin{align*}
\mathbb{P}\left(  \tau^{-}=n\right)   &  =\mathbb{P}\left(  S_{n}\leq
0;\tau^{-}>n-1\right) \\
&  =\int_{0}^{\infty}\mathbb{P}\left(  X_{n}\leq-y\right)  \mathbb{P}\left(
S_{n-1}\in dy;\tau^{-}>n-1\right) \\
&  =\int_{0}^{\varepsilon c_{n}}\mathbb{P}\left(  X\leq-y\right)
\mathbb{P}\left(  S_{n-1}\in dy;\tau^{-}>n-1\right) \\
&  \qquad+\int_{\varepsilon}^{\infty}\mathbb{P}\left(  X\leq-yc_{n}\right)
\mathbb{P}\left(  S_{n-1}\in c_{n}dy;\tau^{-}>n-1\right)  .
\end{align*}
We evaluate the last two integrals separately.

We know from (\ref{Tailtwo}) and (\ref{tailF}) that if $X\in\mathcal{D}\left(
\alpha,\beta\right)  $ with $0<\alpha<2$ and $\beta<1$, then, for a
$q\in(0,1]$,
\begin{equation}
\mathbb{P}\left(  X\leq-y\right)  \sim\frac{q}{y^{\alpha}l_{0}(y)}%
\quad\text{as }y\rightarrow\infty, \label{lefttail}%
\end{equation}
and, according to our construction,%
\[
\mathbb{P}\left(  X\leq-c_{n}\right)  \sim qn^{-1}\quad\text{as }%
n\rightarrow\infty.
\]
Moreover, for any $\varepsilon>0$,
\begin{equation}
\frac{\mathbb{P}\left(  X\leq-yc_{n}\right)  }{\mathbb{P}\left(  X\leq
-c_{n}\right)  }\rightarrow y^{-\alpha}\quad\text{as }n\rightarrow\infty,
\label{aa}%
\end{equation}
uniformly in $y\in(\varepsilon,\infty).$ On the other hand, if $M_{\alpha}%
^{+}(t),0\leq t\leq1,$ is the Levy meander of order $\alpha\neq1$ and the
conditions of Theorem \ref{TlocLad} are valid, then (see \cite{Don85})
\begin{equation}
\left\{  \frac{S_{n}}{c_{n}}\left\vert \tau^{-}>n\right.  \right\}
\overset{d}{\rightarrow}M_{\alpha}^{+}:=M_{\alpha}^{+}(1)\quad\text{as
}n\rightarrow\infty. \label{bb}%
\end{equation}

We show that
\begin{equation}
\int_{0}^{\infty}\frac{\mathbb{P}\left(  M_{\alpha}^{+}\in dy\,\right)
}{y^{\alpha}}<\infty. \label{finInt}%
\end{equation}
Indeed, if this is not the case, for any $N$ one can find $\varepsilon_{N}%
\in(0,1)$ such that%
\[
\int_{\varepsilon_{N}}^{1/\varepsilon_{N}}\frac{\mathbb{P}\left(  M_{\alpha
}^{+}\in dy\,\right)  }{y^{\alpha}}>2N.
\]
This yields%
\begin{align*}
&  \lim_{n\rightarrow\infty}\int_{\varepsilon_{N}}^{1/\varepsilon_{N}}%
\frac{\mathbb{P}\left(  X\leq-yc_{n}\right)  }{\mathbb{P}\left(  X\leq
-c_{n}\right)  }\mathbb{P}\left(  \frac{S_{n-1}}{c_{n}}\in dy\,|\,\tau
^{-}>n-1\right)  \hspace{2cm}\\
&  =\int_{\varepsilon_{N}}^{1/\varepsilon_{N}}\frac{\mathbb{P}\left(
M_{\alpha}^{+}\in dy\,\right)  }{y^{\alpha}}>2N.
\end{align*}
By (\ref{integ1}) we have, as $n\rightarrow\infty,$
\begin{align*}
\frac{2l(n)}{n^{1-\rho}}  &  \geq\mathbb{P}\left(  \tau^{-}>n\right)
=\sum_{k=n+1}^{\infty}\mathbb{P}\left(  \tau^{-}=k\right)  \hspace{7cm}\\
&  \geq\sum_{k=n+1}^{\infty}\mathbb{P}\left(  X_{k}\leq-c_{k}\right)
\mathbb{P}\left(  \tau^{-}>k-1\right)  \times\\
&  \hspace{3cm}\int_{\varepsilon_{N}}^{1/\varepsilon_{N}}\frac{\mathbb{P}%
\left(  X_{k}\leq-yc_{k}\right)  }{\mathbb{P}\left(  X_{k}\leq-c_{k}\right)
}\mathbb{P}\left(  \frac{S_{k-1}}{c_{k}}\in dy\,|\,\tau^{-}>k-1\right) \\
&  \geq N\sum_{k=n+1}^{\infty}\mathbb{P}\left(  X_{k}\leq-c_{k}\right)
\mathbb{P}\left(  \tau^{-}>k-1\right)  \sim N\sum_{k=n+1}^{\infty}\frac
{ql(k)}{k^{2-\rho}}\sim\frac{N}{1-\rho}\frac{ql(n)}{n^{1-\rho}},
\end{align*}
leading to a contradiction for $N>2(1-\rho)q^{-1}$. Thus, (\ref{finInt}) is established.

It easily follows from (\ref{aa}) and (\ref{bb}) that, as $n\rightarrow
\infty,$%
\begin{align}
&  \int_{\varepsilon}^{\infty}\mathbb{P}\left(  X\leq-yc_{n}\right)
\mathbb{P}\left(  S_{n-1}\in c_{n}dy;\tau^{-}>n-1\right) \nonumber\\
&  =\mathbb{P}\left(  X\leq-c_{n}\right)  \mathbb{P}\left(  \tau
^{-}>n-1\right)  \int_{\varepsilon}^{\infty}\frac{\mathbb{P}\left(
X\leq-yc_{n}\right)  }{\mathbb{P}\left(  X\leq-c_{n}\right)  }\mathbb{P}%
\left(  \frac{S_{n-1}}{c_{n}}\in dy\,|\,\tau^{-}>n-1\right) \nonumber\\
&  \sim\frac{ql(n)}{n^{2-\rho}}\int_{\varepsilon}^{\infty}\frac{\mathbb{P}%
\left(  X\leq-yc_{n}\right)  }{\mathbb{P}\left(  X\leq-c_{n}\right)
}\mathbb{P}\left(  \frac{S_{n-1}}{c_{n}}\in dy\,|\,\tau^{-}>n-1\right)
\nonumber\\
&  \sim\frac{ql(n)}{n^{2-\rho}}\int_{\varepsilon}^{\infty}\frac{\mathbb{P}%
\left(  M_{\alpha}^{+}\in dy\,\right)  }{y^{\alpha}}. \label{exact1}%
\end{align}
Taking into account (\ref{finInt}), we obtain
\begin{align}
&  \lim_{\varepsilon\rightarrow0}\lim_{n\rightarrow\infty}\frac{n^{2-\rho}%
}{ql(n)}\int_{\varepsilon}^{\infty}\mathbb{P}\left(  X\leq-yc_{n}\right)
\mathbb{P}\left(  S_{n-1}\in c_{n}dy;\tau^{-}>n-1\right) \nonumber\\
&  =\int_{0}^{\infty}\frac{\mathbb{P}\left(  M_{\alpha}^{+}\in dy\,\right)
}{y^{\alpha}}. \label{exact2}%
\end{align}
To complete the proof of Theorem \ref{TlocLad} it remains to demonstrate that
\begin{equation}
\lim_{\varepsilon\rightarrow0}\limsup_{n\rightarrow\infty}\frac{n^{2-\rho}%
}{l(n)}\int_{0}^{\varepsilon c_{n}}\mathbb{P}\left(  X\leq-y\right)
\mathbb{P}\left(  S_{n-1}\in dy;\tau^{-}>n-1\right)  =0. \label{cc}%
\end{equation}
To this aim we observe that
\begin{align*}
&  \int_{0}^{\varepsilon c_{n}}\mathbb{P}\left(  X\leq-y\right)
\mathbb{P}\left(  S_{n-1}\in dy;\tau^{-}>n-1\right) \\
&  \qquad\qquad\leq\sum_{j=0}^{[\varepsilon c_{n}]+1}\mathbb{P}\left(
X\leq-j\right)  b_{n-1}(j)=:R(\varepsilon c_{n})
\end{align*}
and evaluate $R(\varepsilon c_{n})$ separately for the following three cases:

(i) $0<\alpha<1,$ $|\beta|<1;$

(ii) $1<\alpha<2,$ $|\beta|<1;$

(iii) $1<\alpha<2,$ $\beta=-1.$

\vspace{12pt}

\textbf{(i)}. In view of (\ref{uupb}), (\ref{Tailtwo}) and properties of
regularly varying functions with index $\alpha\in(0,1)$ we have
\begin{align}
R(\varepsilon c_{n})  &  \leq C\frac{1}{c_{n}}\frac{l(n)}{n^{1-\rho}}%
\sum_{j=0}^{[\varepsilon c_{n}]+1}\mathbb{P}\left(  X\leq-j\right) \nonumber\\
&  \leq C_{1}\frac{1}{c_{n}}\frac{l(n)}{n^{1-\rho}}\varepsilon c_{n}%
\mathbb{P}\left(  X\leq-\varepsilon c_{n}\right) \nonumber\\
&  \leq C_{2}\frac{l(n)}{n^{1-\rho}}\varepsilon^{1-\alpha}\frac{l_{0}(c_{n}%
)}{l_{0}(\varepsilon c_{n})}\mathbb{P}\left(  X\leq-c_{n}\right) \nonumber\\
&  \leq C_{3}\frac{l(n)}{n^{2-\rho}}\varepsilon^{1-\alpha}\frac{l_{0}(c_{n}%
)}{l_{0}(\varepsilon c_{n})}\mathbb{\leq}C_{4}\frac{l(n)}{n^{2-\rho}%
}\varepsilon^{1-\alpha-\delta} \label{Remalpha0}%
\end{align}
for any fixed $\delta\in(0,1-\alpha)$ and all sufficiently large $n.$ At the
last step we have used the fact that for every slowly varying function
$l^{\ast}(x)$ and every $\delta>0$ there exists a constant $C_{\delta}$ such
that
\begin{equation}
\frac{l^{\ast}(x)}{l^{\ast}(ax)}\leq C_{\delta}\max\{a^{\delta},a^{-\delta
}\}\quad\text{for all }a,x>0. \label{prop}%
\end{equation}

\textbf{(ii) }In view of (\ref{bsmall2}), equivalences (\ref{RenStand}),
(\ref{Tailtwo}), and estimate (\ref{prop}) with any fixed $\delta\in
(0,\min\{2-\alpha,1-\alpha(1-\rho)\})$, we have for all sufficiently large
$n$,
\begin{align*}
R(\varepsilon c_{n})  &  \leq C\sum_{j=1}^{[\varepsilon c_{n}]+1}\frac
{1}{j^{\alpha}l_{0}(j)}\left(  \frac{j^{\alpha\rho}l_{2}(j)+1}{nc_{n}}%
+\frac{l(n)}{n^{1-\rho}}\frac{j+1}{c_{n}^{2}}\right) \\
&  \leq C_{1}\frac{1}{nc_{n}}\sum_{j=1}^{[\varepsilon c_{n}]+1}\frac{l_{2}%
(j)}{j^{\alpha(1-\rho)}l_{0}(j)}+C\frac{l(n)}{n^{1-\rho}}\frac{1}{c_{n}^{2}%
}\sum_{j=1}^{[\varepsilon c_{n}]+1}\frac{1}{j^{\alpha-1}l_{0}(j)}\\
&  \leq C_{2}\frac{1}{nc_{n}}\left(  \varepsilon c_{n}\right)  ^{1-\alpha
(1-\rho)}\frac{l_{2}(\varepsilon c_{n})}{l_{0}(\varepsilon c_{n})}+C_{3}%
\frac{l(n)}{n^{1-\rho}}\frac{1}{c_{n}^{2}}\frac{\left(  \varepsilon
c_{n}\right)  ^{2-\alpha}}{l_{0}(\varepsilon c_{n})}\\
&  \leq C_{4}\frac{1}{nc_{n}}\left(  \varepsilon c_{n}\right)  ^{1-\alpha
(1-\rho)-\delta}+C_{5}\frac{l(n)}{n^{1-\rho}}\frac{1}{c_{n}^{2}}\left(
\varepsilon c_{n}\right)  ^{2-\alpha-\delta}.
\end{align*}
Hence on account of (\ref{asyma}) we conclude that%
\begin{align}
R(\varepsilon c_{n})  &  \leq C_{4}\frac{\varepsilon^{1-\alpha(1-\rho)-\delta
}}{nc_{n}^{\alpha(1-\rho)+\delta}}+C_{5}\frac{\varepsilon^{2-\alpha-\delta
}l(n)}{n^{1-\rho}}\frac{1}{c_{n}^{\alpha+\delta}}\nonumber\\
&  \leq C_{6}\frac{l(n)}{n^{2-\rho}}\left(  \varepsilon^{1-\alpha
(1-\rho)-\delta}+\varepsilon^{2-\alpha-\delta}\right)  . \label{Remalpha21}%
\end{align}

\textbf{(iii)}. It follows from (\ref{ro}) that if $\beta=-1$, then
$\alpha\rho=1.$ By Lemma \ref{Renew2}, $H(x)\leq Cxl_{3}(x)$. Combining this
estimate with (\ref{bsmall2}), we get
\[
b_{n}(j)\leq C\left(  \frac{jl_{3}(j)+1}{nc_{n}}+\frac{l(n)}{n^{1-\rho}}%
\frac{j+1}{c_{n}^{2}}\right)  .
\]
Recalling (\ref{lefttail}) \ and using (\ref{prop}) once again, we obtain for
any fixed $\delta\in(0,2-\alpha)$ and all $n\geq n(\delta)$,%
\begin{align}
R_{n}(\varepsilon c_{n})  &  \leq C\sum_{j=0}^{[\varepsilon c_{n}%
]+1}\mathbb{P}\left(  X\leq-j\right)  \left(  \frac{jl_{3}(j)+1}{nc_{n}}%
+\frac{l(n)}{n^{1-\rho}}\frac{j+1}{c_{n}^{2}}\right) \nonumber\\
&  \leq C_{1}\left(  \varepsilon c_{n}\right)  ^{2-\alpha}\left(  \frac
{1}{nc_{n}}\frac{l_{3}(\varepsilon c_{n})}{l_{0}(\varepsilon c_{n})}%
+\frac{l(n)}{n^{1-\rho}}\frac{1}{c_{n}^{2}l_{0}(\varepsilon c_{n})}\right)
\nonumber\\
&  \leq C_{2}\varepsilon^{2-\alpha-\delta}\left(  \frac{1}{n}\frac{c_{n}%
l_{3}(c_{n})}{c_{n}^{\alpha}l_{0}(c_{n})}+\frac{l(n)}{n^{1-\rho}}\frac
{1}{c_{n}^{\alpha}l_{0}(c_{n})}\right) \nonumber\\
&  \leq C_{3}\varepsilon^{2-\alpha-\delta}\frac{l(n)}{n^{2-\rho}},
\label{Remalpha3}%
\end{align}
where the inequalities $H(c_{n})\leq Cc_{n}l_{3}(c_{n})\leq Cn^{\rho}l(n)$
have been used for the last step.

Estimates (\ref{Remalpha0}) -- (\ref{Remalpha3}) imply (\ref{cc}). Combining
(\ref{exact2}) with (\ref{cc}) leads to
\begin{equation}
\mathbb{P}\left(  \tau^{-}=n\right)  \sim\frac{ql(n)}{n^{2-\rho}}\int
_{0}^{\infty}\frac{\mathbb{P}\left(  M_{\alpha}^{+}\in dy\,\right)
}{y^{\alpha}}=\frac{ql(n)}{n^{2-\rho}}\mathbb{E}\left(  M_{\alpha}^{+}\right)
^{-\alpha}. \label{local}%
\end{equation}
Summation over $n$ gives%
\[
\mathbb{P}\left(  \tau^{-}>n\right)  =\sum_{k=n+1}^{\infty}\mathbb{P}\left(
\tau^{-}=k\right)  \sim\frac{q}{1-\rho}\frac{l(n)}{n^{1-\rho}}\mathbb{E}%
\left(  M_{\alpha}^{+}\right)  ^{-\alpha}.
\]
Comparing this with (\ref{integ1}), we get an interesting identity%
\begin{equation}
\mathbb{E}\left(  M_{\alpha}^{+}\right)  ^{-\alpha}=(1-\rho)/q
\label{IDENtity}%
\end{equation}
which, in view of (\ref{local}), completes the proof of Theorem \ref{TlocLad}
for $\{0<\alpha<2,\,\beta<1\}\cap\left\{  \alpha\neq1\right\}  $.

\begin{remark}
One can check that the proof of Theorem \ref{TlocLad} for $\{0<\alpha
<2,\,\beta<1\}\cap\left\{  \alpha\neq1\right\}  $ does not use the fact that
in the lattice case the distribution of $X$ is arithmetic.
\end{remark}

%%%%%%%%%%%%%%%%%%%%%%%%%%%%%%%%%%%%%%%%%%%%%%%%%%%%%%%%%%%%%%%%%%%%%%%%%%%%%%%
%%%%%%%%%%%%%%%%%%%%%%%%%%%%%%%%%%%%%%%%%%%%%%%%%%%%%%%%%%%%%%%%%%%%%%%%%%%%%%%

\subsection{Proof of Theorem \ref{TlocLad} for $\{1<\alpha<2,\beta=1\}$}

In view of (\ref{ro}) the assumption $\beta=1$ implies $q=0$ in (\ref{tailF})
and $\rho=1-1/\alpha$. We fix an integer $N>1$ and, for $c_{n}>N,$ write
\begin{align*}
\mathbb{P}\left(  \tau^{-}=n\right)   &  =\int_{0}^{N}\mathbb{P}\left(
X\leq-y\right)  \mathbb{P}\left(  S_{n-1}\in dy;\tau^{-}>n-1\right) \\
&  \qquad+\int_{N}^{c_{n}}\mathbb{P}\left(  X\leq-y\right)  \mathbb{P}\left(
S_{n-1}\in dy;\tau^{-}>n-1\right) \\
&  \qquad\qquad+\int_{c_{n}}^{\infty}\mathbb{P}\left(  X\leq-y\right)
\mathbb{P}\left(  S_{n-1}\in dy;\tau^{-}>n-1\right) \\
&  =:I_{1}(N,n)+I_{2}(N,c_{n})+I_{3}(c_{n}).
\end{align*}
Our aim is to show that the last two integrals divided by $n^{-1/\alpha
-1}l(n)$ vanish as first $n\rightarrow\infty$ and then $N\rightarrow\infty,$
while
\begin{equation}
\lim_{N\rightarrow\infty}\lim_{n\rightarrow\infty}\frac{n^{1+1/\alpha}}%
{l(n)}I_{1}(N,n)=1/\alpha=1-\rho. \label{limInt}%
\end{equation}

To start with, recall that according to Lemma \ref{Lrenewal} under our
conditions
\[
\mathbb{P}(\chi_{+}>x)\sim\frac{1}{x^{\alpha-1}l_{2}(x)}\quad\text{as
}x\rightarrow\infty.
\]
Moreover, it was shown by Doney \cite[Corollary~3]{Don82} that (\ref{IntCond})
is equivalent to the relation $l_{2}(x)\sim Cl_{0}(x)$ as $x\rightarrow\infty
$. Then Lemma \ref{L1} gives the upper bound
\[
b_{n}(x)\leq C\Bigl(\frac{x^{\alpha-1}l_{0}(x)}{nc_{n}}+\frac{l(n)x}%
{n^{1-\rho}c_{n}^{2}}\Bigr)\quad\text{for all }x\geq1.
\]
Besides, Lemma \ref{beta1}, (\ref{asyma}) and (\ref{integ1}) imply exisence of
a constant $K>0$ such that
\begin{equation}
c_{n}\sim\frac{n^{1-\rho}}{Kl(n)}\quad\text{as }n\rightarrow\infty.
\label{A-rel}%
\end{equation}
This equivalence justifies the inequality
\begin{equation}
b_{n}(x)\leq C\frac{l(n)}{n^{2-\rho}}\Bigl(x^{\alpha-1}l_{0}(x)+\frac
{nx}{c_{n}^{2}}\Bigr)\quad\text{for all }x\geq1. \label{B-Bound}%
\end{equation}
As a result, we have for $c_{n}>N>1$ the estimate
\begin{align}
I_{2}(N,c_{n})  &  \leq\sum_{j=N}^{[c_{n}]+1}\mathbb{P}(X\leq-j)b_{n-1}%
(j)\hspace{1cm}\nonumber\\
&  \leq C\frac{l(n)}{n^{2-\rho}}\Bigl(\sum_{j=N}^{[c_{n}]+1}j^{\alpha-1}%
l_{0}(j)\mathbb{P}(X\leq-j)+\frac{n}{c_{n}^{2}}\sum_{j=N}^{[c_{n}%
]+1}j\mathbb{P}(X\leq-j)\Bigr). \label{IntIneq}%
\end{align}
It easily follows from (\ref{IntCond}) and (\ref{tailF}) with $p=1$ and $q=0,$
that
\begin{equation}
\sum_{j=N}^{[c_{n}]+1}j^{\alpha-1}l_{0}(j)\mathbb{P}(X\leq-j)\leq C\sum
_{j=N}^{[c_{n}]+1}\frac{1}{j}\frac{\mathbb{P}(X\leq-j)}{\mathbb{P}(X\geq
j)}\rightarrow0\quad\label{FirstTerm}%
\end{equation}
as first $n\rightarrow\infty$ and than $N\rightarrow\infty.$

Further, recalling that $\mathbb{P}(X\leq-j)=o(\mathbb{P}(X\geq j))$ as
$j\rightarrow\infty,$ we obtain by (\ref{Repa}) and (\ref{tailF}), for
sufficiently large $n$ and a function $r(N)\rightarrow0$ as $N\rightarrow
\infty:$
\begin{align}
\sum_{j=N}^{[c_{n}]+1}j\mathbb{P}(X\leq-j)  &  \leq r(N)\sum_{j=N}^{[c_{n}%
]+1}j\mathbb{P}(X\geq j)\nonumber\\
&  \leq Cr(N)\sum_{j=N}^{[c_{n}]+1}\frac{1}{j^{\alpha-1}l_{0}(j)}\leq
C_{1}r(N)\frac{c_{n}^{2-\alpha}}{l_{0}(c_{n})}\nonumber\\
&  \leq C_{2}r(N)\frac{c_{n}^{2}}{n}. \label{SecondTerm}%
\end{align}
Combining (\ref{IntIneq}), (\ref{FirstTerm}) and (\ref{SecondTerm}), we
conclude that
\begin{equation}
\lim_{N\rightarrow\infty}\limsup_{n\rightarrow\infty}\frac{n^{1+1/\alpha}%
}{l(n)}I_{2}(N,c_{n})=0. \label{Int}%
\end{equation}
To establish a similar result for $I_{3}(c_{n})$, observe that if $\beta=1$,
then, by (\ref{tailF}) and~(\ref{Defa}),
\[
\mathbb{P}(X\leq-c_{n})=o(\mathbb{P}(X\geq c_{n}))=o(1/n)\quad\text{ as
}n\rightarrow\infty,
\]
and, therefore,
\begin{equation}
I_{3}(c_{n})\leq\mathbb{P}(X\leq-c_{n})\mathbb{P}(\tau^{-}>n)=o\Bigl(\frac
{l(n)}{n^{2-\rho}}\Bigr)\quad\text{as }n\rightarrow\infty. \label{Tail}%
\end{equation}

Applying Lemma \ref{L_Conv} and recalling (\ref{A-rel}), we have
\begin{equation}
\lim_{n\rightarrow\infty}\frac{n^{1+1/\alpha}}{l(n)}I_{1}(N,n)=\lim
_{n\rightarrow\infty}Knc_{n}I_{1}(N,n)=K\int_{0}^{N}\mathbf{P}(X\leq
-x)\mu(dx). \label{IntN}%
\end{equation}
In view of (\ref{B-Bound}),
\[
\mu((x,x+1])=\lim_{n\rightarrow\infty}nc_{n}b_{n}(x)\leq Cx^{\alpha-1}\ell
_{0}(x).
\]
From this, taking into account conditions (\ref{FirstTerm}) and (\ref{IntCond}%
), we get%
\[
\int_{0}^{\infty}\mathbf{P}(X\leq-x)\mu(dx)<\infty.
\]
Hence we conclude that
\begin{equation}
\lim_{N\rightarrow\infty}\lim_{n\rightarrow\infty}\frac{n^{1+1/\alpha}}%
{l(n)}I_{1}(N,n)=K\int_{0}^{\infty}\mathbf{P}(X\leq-x)\mu(dx).
\label{MainTerm}%
\end{equation}
Combining (\ref{Int}), (\ref{Tail}) and (\ref{MainTerm}) yields, as
$n\rightarrow\infty$,
\begin{equation}
\mathbb{P}(\tau^{-}=n)\sim\frac{Kl(n)}{n^{1+1/\alpha}}\int_{0}^{\infty
}\mathbf{P}(X\leq-x)\mu(dx)\sim\frac{1}{nc_{n}}\int_{0}^{\infty}%
\mathbf{P}(X\leq-x)\mu(dx)\text{ }. \label{Intmu}%
\end{equation}
Comparing this formula with the tail behavior of $\tau^{-}$ given by
(\ref{integ1}) leads to the equalities
\begin{equation}
K\int_{0}^{\infty}\mathbf{P}(X\leq-x)\mu(dx)=1-\rho=1/\alpha. \label{INtalpha}%
\end{equation}
This justifies (\ref{limInt}), finishing the proof of our theorem for
$1<\alpha<2,\beta=1$.
%%%%%%%%%%%%%%%%%%%%%%%%%%%%%%%%%%%%%%%%%%%%%%%%%%%%%%%%%%%%%%%%%%%%%%%%%%%%%%%%%%%%%%%%%%%
%%%%%%%%%%%%%%%%%%%%%%%%%%%%%%%%%%%%%%%%%%%%%%%%%%%%%%%%%%%%%%%%%%%%%%%%%%%%%%%%%%%%%%%%%%%

\subsection{Proof of Theorem \ref{TlocLad} for $\left\{  \alpha=2,\beta
=0\right\}  $}

Consider first the case of arithmetic distributions and assume for simplicity
that $h=1$ from now on. In this case we write%
\begin{align*}
\mathbb{P}(\tau^{-}  &  =n)=\sum_{j=1}^{\infty}\mathbb{P}(X\leq-j)\mathbb{P}%
(S_{n-1}=j;\tau^{-}>n-1)\\
&  =\Delta_{1}(c_{n})+\Delta_{2}(c_{n}),
\end{align*}
where%
\[
\Delta_{1}(c_{n}):=\sum_{j=1}^{\left[  c_{n}\right]  }\mathbb{P}%
(X\leq-j)\mathbb{P}(S_{n-1}=j;\tau^{-}>n-1),
\]%
\[
\Delta_{2}(c_{n}):=\sum_{j=\left[  c_{n}\right]  +1}^{\infty}\mathbb{P}%
(X\leq-j)\mathbb{P}(S_{n-1}=j;\tau^{-}>n-1).
\]
Recall that if $\alpha=2$ then $\rho=1/2$. In view of (\ref{102}),
(\ref{ProV}) and (\ref{integ1})
\begin{align*}
\Delta_{2}(c_{n})  &  \leq\mathbb{P}(X\leq-c_{n})\mathbb{P}(\tau^{-}>n-1)\\
&  =o\left(  \frac{1}{n}\frac{l(n)}{n^{1/2}}\right)  =\text{ }o\left(
\frac{l(n)}{n^{3/2}}\right)  \text{\ as \ }n\rightarrow\infty.
\end{align*}
To evaluate $\Delta_{1}(c_{n})$ denote $g_{2,0}(x)=(\sqrt{2\pi})^{-1}%
\exp\left\{  -x^{2}/2\right\}  ,\ x\in(-\infty,\infty),$ the density of the
standard normal law and set%
\[
w(n):=\sum_{j=1}^{\left[  c_{n}\right]  }g_{2,0}\left(  \frac{j}{c_{n}%
}\right)  \mathbb{P}(X\leq-j)H(j-1).
\]
By formula (3.15) in \cite{BJD06}, as $n\rightarrow\infty$,%
\[
\mathbb{P}(S_{n-1}=j;\tau^{-}>n-1)\sim\frac{H(j-1)}{n}\mathbb{P}%
(S_{n-1}=j)\sim\frac{H(j-1)}{nc_{n}}g_{2,0}\left(  \frac{j}{c_{n}}\right)
\]
uniformly in $j\in\lbrack1,c_{n}].$ This gives
\begin{equation}
\Delta_{1}(c_{n})=\frac{1+r(n)}{nc_{n}}w(n), \label{65}%
\end{equation}
where $\ r(n)\rightarrow0$ as $n\rightarrow\infty.$ As a result we obtain%
\begin{equation}
\mathbb{P}(\tau^{-}=n)=\frac{1+r(n)}{nc_{n}}w(n)+o\left(  \frac{l(n)}{n^{3/2}%
}\right)  . \label{Prelim}%
\end{equation}
Hence it follows that, as $n\rightarrow\infty$,%
\begin{align*}
\frac{l(n)}{n^{1/2}}  &  \sim\mathbb{P}(\tau^{-}>n)=\sum_{k=n+1}^{\infty
}\left(  \frac{1+r(k)}{kc_{k}}w(k)+o\left(  \frac{l(k)}{k^{3/2}}\right)
\right) \\
&  =\left(  1+r_{1}(n)\right)  \sum_{k=n+1}^{\infty}\frac{w(k)}{kc_{k}%
}+o\left(  \frac{l(n)}{n^{1/2}}\right)  ,
\end{align*}
where $r_{1}(n)\rightarrow0$ as $n\rightarrow\infty.$ Since $w(n)$ is monotone
increasing in $n$, and $c_{n}\sim n^{1/2}l_{1}(n)$ as $n\rightarrow\infty,$
Lemma \ref{Ltec} with $\gamma=1-\rho=1/2$ yields after obvious
transformations
\begin{equation}
\frac{w(n)}{nc_{n}}\sim\frac{1}{2}\frac{l(n)}{n^{3/2}}\text{ \ as
\ }n\rightarrow\infty, \label{67}%
\end{equation}
which, on account of (\ref{Prelim}) \ finishes the proof of (\ref{main}) for
$\left\{  \alpha=2,\beta=0\right\}  $ in the arithmetic case. To establish the
same result for non-lattice distributions one should apply the respective
statements in \cite{BJ04}.
%\newpage

%%%%%%%%%%%%%%%%%%%%%%%%%%%%%%%%%%%%%%%%%%%%%%%%%%%%%%%%%%%%%%%%%%%%%%%%%%%%%%%%%
%%%%%%%%%%%%%%%%%%%%%%%%%%%%%%%%%%%%%%%%%%%%%%%%%%%%%%%%%%%%%%%%%%%%%%%%%%%%%%%%%

\section{Proof of Theorem \ref{propos}}

Applying (\ref{strict}) to the random walk $\{-S_{n}\}_{n\geq0}$, we have
\[
1-\mathbb{E}z^{T^{-}}=\exp\left\{  -\sum_{n=1}^{\infty}\frac{z^{n}}%
{n}\mathbb{P}(S_{n}<0)\right\}  .
\]
Recalling (\ref{defOmega}) and (\ref{weak}) we obtain
\begin{equation}
1-\mathbb{E}z^{T^{-}}=\left(  1-\mathbb{E}z^{\tau^{-}}\right)  \Omega(z).
\label{Rec}%
\end{equation}
On account of $\mathbb{P}(\tau^{-}=0)=0$, equality (\ref{Rec}) implies
\begin{equation}
\mathbb{P}(T^{-}=n)=\sum_{k=1}^{n}\mathbb{P}(\tau^{-}=k)\omega_{n-k}%
-\omega_{n},\ n\geq1. \label{Recur1}%
\end{equation}
Suppose first that the distribution of $X$ is arithmetic. By the Gnedenko
local theorem we get for this case
\[
\frac{1}{n}\mathbb{P}(S_{n}=0)=\frac{g_{\alpha,\beta}(0)}{nc_{n}}%
(1+o(1))\quad\text{as }n\rightarrow\infty.
\]
This representation and Theorem 2 in \cite{CNW} \ provide existence of a
constant $C>0$ such that
\[
\omega_{n}=\frac{C}{nc_{n}}(1+o(1))\quad\text{as }n\rightarrow\infty.
\]
Using this equality and (\ref{main}) in (\ref{Recur1}) and recalling that
$\mathbb{P}(\tau^{-}<\infty)=1$, we obtain
\[
\mathbb{P}(T^{-}=n)=\Omega(1)\mathbb{P}(\tau^{-}=n)(1+o(1))+o((nc_{n}%
)^{-1})\quad\text{as }n\rightarrow\infty.
\]
Observing that $\mathbb{P}(\tau^{-}=n)\geq C/nc_{n}$, we get the desired
statement for the arithmetic case.

If the distribution of $X$ is non-lattice, then there exists a constant
$r\in(0,1)$ such that $\mathbb{P}(S_{n}=0)\leq r^{n}$ for all $n\geq1$ (we may
choose $r$ as the total mass of the lattice component of the distribution of
$X$). Consequently, $\omega_{n}\leq r^{n}$ for all $n\geq1$. From this
estimate and (\ref{Recur1}) we see that the statement of Theorem \ref{propos}
is valid in the non-lattice case as well.

\section{Discussion and concluding remarks}

We see by (\ref{weak}) that the distribution of \ $\tau^{-}$ is completely
specified by the sequence $\{\mathbb{P}\left(  \,S_{n}>0\right)  \}_{n\geq1}$.
As we have mentioned in the introduction, the validity of condition
(\ref{SpitDon}) is sufficient to reveal the asymptotic behavior of
$\mathbb{P}(\tau^{-}>n)$ as $n\rightarrow\infty$. Thus, \ in view of
(\ref{integ1}), nonformal arguments based on the plausible smoothness of
$l(n)$ immediately give the desired answer%
\begin{align*}
\mathbb{P}(\tau^{-}  &  =n)=\mathbb{P}(\tau^{-}>n-1)-\mathbb{P}(\tau^{-}>n)\\
&  =\frac{l(n-1)}{\left(  n-1\right)  ^{1-\rho}}-\frac{l(n)}{n^{1-\rho}%
}\approx l(n)\left(  \frac{1}{\left(  n-1\right)  ^{1-\rho}}-\frac
{1}{n^{1-\rho}}\right) \\
&  \approx\frac{(1-\rho)l(n)}{n^{2-\rho}}\sim\frac{1-\rho}{n}\mathbb{P}%
(\tau^{-}>n)
\end{align*}
under the Doney condition only. In the present paper we failed to achieve such
a generality. However, it is worth to be mentioned that the Doney condition,
being formally weaker than the conditions of Theorem \ref{TlocLad}, requires
in the general case the knowledge of the behavior of the whole sequence
$\{\mathbb{P}\left(  \,S_{n}>0\right)  \}_{n\geq1},$ while the assumptions of
Theorem \ref{TlocLad} concern a single summand only. \ Of course, imposing a
stronger condition makes our life easier and allows us to give, in a sense, a
constructive proof showing what happens in reality at the distant moment
$\tau^{-}$ of the first jump of the random walk in question below zero.
Indeed, our arguments for the case $\left\{  0<\alpha<2,\ \beta<1\right\}
\cap\left\{  \alpha\neq1\right\}  $ demonstrate (compare (\ref{lefttail}),
(\ref{aa}), and (\ref{exact1})) that for any $x_{2}>x_{1}>0$,
\begin{align*}
\lim_{n\rightarrow\infty}\mathbb{P}(S_{n-1}  &  \in(c_{n}x_{1},c_{n}%
x_{2}]|\tau^{-}=n)\\
&  \hspace{-1cm} =\lim_{n\rightarrow\infty}\frac{\mathbb{P}(\tau^{-}%
>n-1)}{\mathbb{P}(\tau^{-}=n)}\int_{x_{1}}^{x_{2}}\mathbb{P}(X<-yc_{n}%
)\mathbb{P}(S_{n-1}\in c_{n}dy|\tau^{-}>n-1)\\
&  \hspace{-1cm} =\lim_{n\rightarrow\infty}\frac{\mathbb{P}(\tau^{-}%
>n-1)q}{\mathbb{P}(\tau^{-}=n)n}\int_{x_{1}}^{x_{2}}\frac{\mathbb{P}%
(X<-yc_{n})}{\mathbb{P}(X<-c_{n})}\mathbb{P}(S_{n-1}\in c_{n}dy|\tau
^{-}>n-1)\\
&  \hspace{-1cm} =\frac{q}{1-\rho}\int_{x_{1}}^{x_{2}}\frac{\mathbb{P}%
(M_{\alpha}^{+}\in dy)}{y^{\alpha}}.
\end{align*}
In view of (\ref{IDENtity}) this means that the contribution of the
trajectories of the random walk satisfying $S_{n-1}c_{n}^{-1}\rightarrow0$ or
$S_{n-1}c_{n}^{-1}\rightarrow\infty$ as $n\rightarrow\infty$ to the event
$\left\{  \tau^{-}=n\right\}  $ is negligibly small in probability. A
"typical" trajectory looks in this case as follows: it is located over the
level zero up to moment $n-1$ with $S_{n-1}\in(\varepsilon c_{n}$
$,\varepsilon^{-1}c_{n})$ for sufficiently small $\varepsilon>0$ and at moment
$\tau^{-}=n$ the trajectory makes a big negative jump $X_{n}<-S_{n-1}$ of
order $O(c_{n}).$

On the other hand, if $\left\{  1<\alpha<2,\ \beta=1\right\}  $ and condition
(\ref{IntCond}) holds, then (compare (\ref{Conv}), (\ref{IntN}),
(\ref{Intmu}), and (\ref{INtalpha})) for any $N_{2}>N_{1}>0$,
\begin{align*}
\lim_{n\rightarrow\infty}\mathbb{P}(S_{n-1}  &  \in(N_{1},N_{2}]|\tau^{-}=n)\\
&  =\lim_{n\rightarrow\infty}\frac{1}{\mathbb{P}(\tau^{-}=n)}\int_{N_{1}%
}^{N_{2}}\mathbb{P}(X<-y)\mathbb{P}(S_{n-1}\in dy;\tau^{-}>n-1)\\
&  =\lim_{n\rightarrow\infty}K\alpha nc_{n}\int_{N_{1}}^{N_{2}}\mathbb{P}%
(X<-y)\mathbb{P}(S_{n-1}\in dy;\tau^{-}>n-1)\\
&  =K\alpha\int_{N_{1}}^{N_{2}}\mathbb{P}(X<-y)\mu(dy).
\end{align*}
Thus, the main contribution to $\mathbb{P}\left(  \tau^{-}=n\right)  $ is
given in this case by the trajectories located over the level zero up to
moment $n-1$ with $S_{n-1}\in\lbrack0,N]$ for sufficiently big $N$ and with
not "too big" jump $X_{n}<-S_{n-1}$ of order $O(1).$

Unfortunately, our approach to investigate the behavior of $\mathbb{P}%
(\tau^{-}=n)$ in the case $\alpha=2$ is pure analytical and does not allow us
to extract typical trajectories without further restrictions on the
distribution of $X$. However, we can still deduce from our proof some
properties of the random walk conditioned on $\{\tau^{-}=n\}$. Observe that,
for any fixed $\varepsilon>0$, the trajectories with $S_{n-1}>\varepsilon
c_{n}$ give no essential contribution to $\mathbb{P}(\tau^{-}=n)$. Indeed, it
follows from (\ref{65}) and (\ref{67}) that $\Delta_{1}(\varepsilon c_{n}%
)\sim\Delta_{1}(c_{n})$ as $n\rightarrow\infty$ for every fixed $\varepsilon$.
This, along with the estimate from above for $\Delta_{2}(c_{n}),$ gives the
claimed property. Furthermore, one can easily verify that if $\sum
_{j=1}^{\infty}\mathbb{P}(X\leq-j)H(j)=\infty$, then for every $N\geq1$,
\[
\sum_{j=1}^{N}\mathbb{P}(X\leq-j)\mathbb{P}(S_{n-1}=j;\tau^{-}>n-1)=o\left(
\frac{l(n)}{n^{3/2}}\right)  \quad\text{as }n\rightarrow\infty,
\]
i.e. the contribution of the trajectories with $S_{n-1}=O(1)$ to
$\mathbb{P}(\tau^{-}=n)$ is negligible small. As a result we see that
$S_{n-1}\rightarrow\infty$ but $S_{n-1}=o(c_{n})$ for all "typical"
trajectories meeting the condition $\{\tau^{-}=n\}$. Thus, under the
conditions of Theorem \ref{TlocLad} we have for $\alpha=2$ a kind of
"continuous transition" between the two strategies that take place for the
case $\alpha<2$. We note, for completeness, that if $\sum_{j=1}^{\infty
}\mathbb{P}(X\leq-j)H(j)$ is finite, then the typical behavior of the
trajectories is similar to that for the case $\left\{  0<\alpha<2,\ \beta
=1\right\}  $.

Unfortunately, the methods of the present paper do not work for $\alpha=1$,
and we leave the problem on the asymptotic behavior of $\mathbb{P}(\tau
^{-}=n)$ open for this case.

\vspace*{12pt} \textit{Acknowledgement}. The main results of the paper were
obtained during visits of the first author to the Weierstrass Institute in
Berlin and the second author to the Steklov Mathematical Institute in Moscow.
The hospitality of the both institutes is greatly acknowledged.

%%%%%%%%%%%%%%%%%%%%%%%%%%%%%%%%%%%%%%%%%%%%%%%%%%%%%%%%%%%%%%%%%%%%%%%%%%%%%%%%

\end{document}